\documentclass{article}
\usepackage{amssymb}
\usepackage{amsfonts}
\usepackage{amsmath}
\usepackage{makeidx}

\newtheorem{theorem}{Theorem}

\newtheorem{corollary}[theorem]{Corollary}

\newtheorem{definition}[theorem]{Definition}
\newtheorem{example}[theorem]{Example}

\newtheorem{proposition}[theorem]{Proposition}
\newtheorem{remark}[theorem]{Remark}

\newenvironment{proof}[1][Proof]{\noindent\textbf{#1.} }{\ \rule{0.5em}{0.5em}}

\begin{document}

\begin{center}
{\Large The lattice of ideals of certain rings }

\begin{equation*}
\end{equation*}%
Diana SAVIN \\

Transilvania University of Bra\c {s}ov,\\ Faculty of Mathematics and Computer Science, 
 Bra\c {s}ov, Romania\\
e-mail: diana.savin@unitbv.ro

\begin{equation*}
\end{equation*}
\end{center}

\textbf{Abstract. }{\small Let $A$ be a unitary ring and  let $(\mathbf{I(A),\subseteq })$  be the lattice of ideals of the ring $A.$ In this article we will study the property of the lattice $(\mathbf{I(A),\subseteq })$ to be Noetherian or not, for various types of rings $A$. 
 In the last section of the article we study certain rings that are not Boolean rings, but all their ideals are idempotent.}%

\begin{equation*}
\end{equation*}


\medskip


\section{Preliminaries}

Let  $\left(A, +, \cdot\right)$ be a unitary ring and let $\mathbf{I(A)}$ be the set of ideals of the ring $A.$ It is known that $(\mathbf{I(A),\subseteq })$  is a complete modular lattice.\\
 Lattice theory has connections and applications not only in mathematical logic but also in other branches of algebra such as ring theory, group theory, the theory of algebraic hyperstructures and algebraic number theory (see \cite{ps}, \cite{Ciungu}, \cite{Ighedo}, \cite{piciu}).
Also, lattice theory has applications in discrete mathematics, applied mathematics, physics and computer science (see \cite{af}, \cite{Nuland}).

Many authors studied the lattice of ideals of some ring, determining many of its properties (see \cite{Blair}, \cite{Ighedo}).
In this paper we study some properties of ideals lattice of some special rings. The structure of the article is as follows: the first section is a preliminary results section. In the second section we study the lattice of ideals of a Noetherian ring, of a principal ring, of a Bezout ring, of a Dedekind ring. In the 3rd section we study the lattice of ideals of the ring of algebraic integers of a quadratic field, in case when this ring is not a principal ring.
In the last section we determine some interesting properties of the lattice of ideals for rings of residue classes modulo $n$ (whith $n$ a positive integer, $n$$\geq 2$).\\

We recall some results about lattices, about Bezout rings and about ring of algebraic integers, which
will be necessary in our proofs.

\begin{definition} (\cite{nas}). 
\label{*}
A partial ordered set $(A,\leq )$ is Noetherian if every increasing chain $%
x_{1}\leq x_{2}\leq ...$ of elements of $A$ is stationary, i.e., there is a
positive integer $n\geq 1$ such that $x_{n+i}=x_{n},$ for every $i\geq 0$. 
\end{definition}
\begin{proposition}
\label{prop21.} (\cite{al}, \cite{sa}).  Let  $\left(A, +, \cdot\right)$ be a principal ring. Then:
\begin{enumerate}
\item[$i)$] $a\cdot A+ b\cdot A=\left(a,b\right)\cdot A;$

\item[$ii)$] $a\cdot A \cap b\cdot A=\left[a,b\right]\cdot A.$

\end{enumerate}
\end{proposition}
We recall the following definition.

\begin{definition} (\cite{al}). 
\label{Bezout1}
A Bezout domain is an integral domain in which every finitely generated ideal is principal.
\end{definition}
A result similar to the one in Proposition \ref{prop21.} also occurs in Bezout rings, but only for finitely generated ideals.
\begin{proposition}  (\cite{co}).
\label{prop22.}  Let  $\left(A, +, \cdot\right)$ be a Bezout ring. Then:
\begin{enumerate}
\item[$i)$] $a\cdot A+ b\cdot A=\left(a,b\right)\cdot A;$

\item[$ii)$] $a\cdot A \cap b\cdot A=\left[a,b\right]\cdot A,$
\end{enumerate}
for any two finitely generated ideals $I$ and $J$  of the ring $A$ (so $I=a\cdot A$ and $J=b\cdot A,$ 
since $A$ is a Bezout ring).
\end{proposition}

Next, we recall some known results about quadratic fields, results that we will use in our proofs.\\
Let a quadratic field  $K=\mathbb{Q}(\sqrt{d}),$ where $d\neq 0,1$  is a square free   integer. We recall that the discriminant of $K$ is $$\displaystyle\Delta_{K}=\left\{
\begin{array}{rl}
4d, & \hbox{if $d$$\equiv$ $2$ or $3$ (mod $4$);} \\
d, & \hbox{if $d$$\equiv$ $1$ (mod $4$).}
\end{array}
\right.
$$
(see  \cite{al}, \cite{ireland},  \cite{sa})\\

The following theorem gives us the decomposition of a prime integer in the ring of algebraic integers of a quadratic field.

\begin{theorem}  (\cite{Ribenboim}, \cite{al}, \cite{ireland},  \cite{sa}).
\label{theorem23.}
Let $d\neq 0,1$   {be a square free   integer.}   {Let} $\Delta_{K}$ {be the discriminant of} $K$ {and let} $\mathcal{O}_{K}$   {be the ring of integers of the quadratic field} $K=\mathbb{Q}(\sqrt{d})$. {Let} $p$  {be a prime integer}, $p\geq 3.$ {Then, we have}:
\begin{enumerate}
\item[$i)$] $p$   {is ramified in} $\mathcal{O}_{K}$ if and only if
$p  | \Delta_{K}$. In this case 
$p\mathcal{O}_{K}=(p   , \sqrt{d}   )^2$;
\item[$ii)$]  $p$   splits totally   in $\mathcal{O}_{K}$  if and only if
  $(\frac{\Delta_{K}}{p})=1$. In this case $p\mathcal{O}_{K}=P \cdot \overline{P},$   {where} $P,$ $\overline{P}$ $\in$ Spec$\left(\mathcal{O}_{K}\right),$ $P=\left(p, a+\sqrt{d}\right),$  $\overline{P}=\left(p, a-\sqrt{d}\right),$ where $a$$\in$$\mathbb{Z}$ is a solution of the congruence $x^{2}$$\equiv$$d$ (mod $p$);
\item[$iii)$]  $p$   {is inert in} $\mathcal{O}_{K}$ {if and only if}
 $(\frac{\Delta_{K}}{p})=-1$; 
\item[$iv)$] 
 2 is ramified in $\mathcal{O}_{K}$  {if and only if}
 $d\equiv 2$ (mod $4$) or $d\equiv 3$ (mod $4$) 
  In the first case    $2\mathcal{O}_{K}=(2,\sqrt{d} )^{2}$, while in the second case
 $2\mathcal{O}_{K}=(2,1+\sqrt{d} )^{2};$
\item[$v)$] 
$2$  splits totally   in $\mathcal{O}_{K}$  {if and only if} 
$d\equiv 1$ (mod 8). In this case     $2\mathcal{O}_{K}=P\cdot P^{'}$,   {where} $P, P^{'}$$\in$$\mathcal{O}_{K}$, 
$P\neq P^{'},$ $P=(2, \frac{1+\sqrt{d}}{2})$;
\item[$vi)$] 
the prime $2 $   {is inert  in} $\mathcal{O}_{K}$  {if and only if} 
$d\equiv 5$ (mod $8$). 
\end{enumerate}
\end{theorem}

\begin{equation*}
\end{equation*}%
\section{The lattice of ideals of principal rings, Noetherian rings, Bezout rings, Dedekind rings }

\begin{equation*}
\end{equation*}%
\begin{remark}
\label{remark31.}  Let  $\left(A, +, \cdot\right)$ be a principal ring. Then, the lattice $(\mathbf{I(A),\subseteq })$ of the ideals of the ring $A$ is a Noetherian lattice.
\end{remark}

\begin{proof} 
Since $A$ is a principal ring, it results that $A$ is a Noetherian ring, so the lattice $(\mathbf{I(A),\subseteq })$ of the ideals of the ring $A$ is a Noetherian lattice.
\end{proof}

We immediately obtain the following corollary.
\begin{corollary}
\label{corollary32.}  The lattice $(\mathbf{I(\mathbb{Z}),\subseteq })$ of the ideals of the ring of the integer numbers is Noetherian.
\end{corollary}

\begin{proposition}
\label{prop33.}  Let  $\left(A, +, \cdot\right)$ be a Bezout ring and let $B$ the set of all finitely generated ideals of the ring $A.$ Then:
\begin{enumerate}
\item[$i)$] $B$ is a sublattice of the the lattice $(\mathbf{I(A),\subseteq })$ of the ideals of the ring $A;$ 
\item[$ii)$] $B$ is a Noetherian lattice.
\end{enumerate}
\end{proposition}
\begin{proof}
 i) Let $I$ and $J$ be two  finitely generated ideals of the ring $A.$ Since $A$ is a Bezout ring, it results that there are $a, b$$\in$$A$ such that $I=a\cdot A$ and  $J=b\cdot A.$ Applying  Proposition \ref{prop22.} we ovtain that 
$I+J=\left(a,b\right)\cdot A$$\in$$B$ and  $I \cap  J=\left[a,b\right]\cdot A$$\in$$B,$ so, $B$ is a sublattice of the the lattice $(\mathbf{I(A),\subseteq }).$\\
ii) Let $I$ be an ideal of the ring $\left(B, +, \cdot\right),$ so it is a principal ideal of the ring $A.$\\
 If $I=p^{n}\cdot A,$ where $n$$\in$$\mathbb{N}^{*}$ and $p$ is a prime element of the ring  $\left(A, +, \cdot\right)$, then the following increasing chain of ideals of the ring  $\left(A, +, \cdot\right)$:
$$ p^{n}\cdot A\subseteq p^{n-1}\cdot A\subseteq...\subseteq p\cdot A  $$
is stationary.\\
Let  $U\left( A \right)$ be the set of invertible elements of the ring $A.$\\
If $I=a\cdot A,$ where $a$$\in$ $A\backslash U\left( A \right),$ $a \neq 0.$ Let $p_{1}$ be a prime element from the ring $A,$ with the property $p_{1} | A.$  
Let $\alpha_{1}$ be the greatest positive integer with the property $ p^{\alpha_{1}}_{1}| A.$  Then the following increasing chain of ideals of the ring  $\left(A, +, \cdot\right)$:
$$a\cdot A\subseteq p^{\alpha_{1}}_{1}\cdot A\subseteq   p^{\alpha_{1}-1}_{1}\cdot A\subseteq...\subseteq p_{1}\cdot A$$
is stationary.
\end{proof}
\smallskip\\
There are Bezout rings that also contain infinitely generated ideals, i.e. generally Bezout rings that are not Noetherian rings. I now give such an example.
\begin{remark}  (\cite{al}) 
\label{remark34.}  The ring $\left(\mathcal{O}_{\mathbb{C}}, +, \cdot\right)$  of all algebraic integers is a Bezout ring, but it is not a Noetherian ring, so it is not a Dedekind ring.
\end{remark}

\begin{corollary}
\label{corollary35.}  The lattice $(\mathbf{I(\mathcal{O}_{\mathbb{C}}),\subseteq })$ of the ideals of the ring of all algebraic integers is not Noetherian.
\end{corollary}
\begin{proof} It follows immediately from the Remark  \ref{remark34.}. \\
We now give a direct proof of this remark.
Let $n$ be a positive integer, $n\geq 2.$   We remark that $\sqrt[2^{n}]{5}$ is a root of the polynomial $f=X^{2^n} - 5,$ so  $\sqrt[2^{n}]{5}$$\in$$\mathcal{O}_{\mathbb{C}}.$\\
In the ring $\left(\mathcal{O}_{\mathbb{C}}, +, \cdot\right)$ we consider the following  increasing chain of ideals:
$$ \sqrt{5}\cdot \mathcal{O}_{\mathbb{C}}\subseteq \sqrt[4]{5}\cdot \mathcal{O}_{\mathbb{C}}\subseteq...\subseteq \sqrt[2^{n}]{5}\cdot \mathcal{O}_{\mathbb{C}}.\subseteq...  $$
which is not stationary. It resullts that the lattice $(\mathbf{I(\mathcal{O}_{\mathbb{C}}),\subseteq })$ is not Noetherian.
\end{proof}

From those  previously proved, the following corollary immediately follows.
\begin{corollary} \label{corollary36.}
The following statements are true:
\begin{enumerate}
\item[$i)$] The lattice of ideals of any principal ring is a Noetherian lattice;
\item[$ii)$] The lattice of ideals of any Dedekind ring is a Noetherian lattice;
\item[$iii)$] There are Bezout rings whose lattice of ideals is not a Noetherian lattice;
\item[$iv)$] The lattice of ideals of any finite Bezout ring is a Noetherian lattice.
\end{enumerate}
\end{corollary}

\section{Rings with the lattice of ideals being Noetherian, with the property that there are non-principal ideals whose product is a principal ideal}

In \cite{ps} (Theorem 2), D. Piciu and D. Savin obtained that:\\
 \textit {A residuated lattice} L  \textit {is i-Noetherian (this means that the ideals lattice} $(\mathbf{I(L),\subseteq })$ \textit {is Noetherian)  if and only if every ideal of} \textit {L is principal}.\\
This result is not true if instead of a residuated lattice $L$, we consider a unitary ring $A.$ It is easy to find such rings. In the following we find some examples of non-principal rings, with the lattice of ideals to be Noetherian, but in which the product of certain non-principal ideals is a principal ideal.
We found such examples among rings of algebraic integers.\\
Let $A$ be a unitary commutative ring and let $I$ be an ideal of the ring $A$ gene-\\
rated by the elements $\alpha_{1},$  $\alpha_{2},$ ....,  $\alpha_{n}$$\in$$A,$ where $n$ a positive integer, $n\geq 2$. We denote the ideal $I$ as follows: $I=\left(\alpha_{1}, \alpha_{2},\ldots,\alpha_{n}\right).$  \\
We know that the lattice of ideals of a ring of algebraic integers (of an algebraic numbers field $K$) is a distributive lattice (see \cite{Blair}).\\
We start from the following example (from (\cite{al})): let the quadratic field $\mathbb{Q}\left(i\cdot \sqrt{5}\right)$ (where $i^{2}=-1$) and its ring of algebraic integers $\mathbb{Z}\left[i\cdot \sqrt{5}\right]$. In this ring we consider the prime ideals $I_{1}=\left(3, 4+i\cdot \sqrt{5}\right),$  $I_{2}=\left(3, 4- i\cdot \sqrt{5}\right),$
$J_{1}=\left(7, 4+i\cdot \sqrt{5}\right),$  $J_{2}=\left(7, 4- i\cdot \sqrt{5}\right).$ Atunci $I_{1}\cdot I_{2}=\left(3\right), $ $J_{1}\cdot J_{2}=\left(7\right), $ $I_{1}\cdot J_{1}=\left(4+i\cdot \sqrt{5}\right), $ $I_{2}\cdot J_{2}=\left(4-i\cdot \sqrt{5}\right).$
So, $I_{1}, $ $I_{2}, $  $J_{1}, $ $J_{2} $ are non-principal  ideals of the ring $\mathbb{Z}\left[i\cdot \sqrt{5}\right],$ but  $I_{1}\cdot I_{2},$  $I_{1}\cdot J_{1},$  $I_{2}\cdot J_{2},$ $J_{1}\cdot J_{2}$ are principal ideals of the ring $\mathbb{Z}\left[i\cdot \sqrt{5}\right].$\\
We obtain the following result:
\begin{theorem}  
\label{theorem37.}
Let $d\neq 0,1$   {be a square free  integer and let a quadratic field}   $K=\mathbb{Q}(\sqrt{d}).$   {Let} $\Delta_{K}$ {be the discriminant of} $K$ {and let} $\mathcal{O}_{K}$   {be the ring of integers of the field} $K$ such that  $\mathcal{O}_{K}$  is not principal ring.  {Then, we have}:
\begin{enumerate}
\item[$i)$]   {if} $p$  {is an odd prime integer},  then the ideal $p\mathcal{O}_{K}=P_1 \cdot P_2,$ where $P_1, P_2$$\in$ Spec$\left(\mathcal{O}_{K}\right)$  are non-principal ideals (each of the ideals $P_1$ and $P_2$ has $2$ generators) if and only if $d$ is a quadratic residue modulo $p$ or $p  | d.$
\item[$ii)$]  {if} $p$  {is an odd prime integer},  then the ideal $p\mathcal{O}_{K}=P_1 \cdot P_2,$ where $P_1, P_2$$\in$ Spec$\left(\mathcal{O}_{K}\right),$   $P_1 \neq P_2$  are non-principal ideals (each of the ideals $P_1$ and $P_2$ has $2$ generators) if and only if $d$ is a quadratic residue modulo $p.$ 
In this case $p\mathcal{O}_{K}=P_1 \cdot P_2=P_1 \cap P_2.$ 
\item[$iii)$]  the ideal  $2\mathcal{O}_{K}=P_1 \cdot P_2,$ where $P_1, P_2$$\in$ Spec$\left(\mathcal{O}_{K}\right)$ are non-principal ideals (each of the ideals $P_1$ and $P_2$ has $2$ generators) if and only if $d\equiv 2$ (mod $4$) or $d\equiv 3$ (mod $4$) or $d\equiv 1$ (mod $8$) 
\item[$iv)$]  the ideal  $2\mathcal{O}_{K}=P_1 \cdot P_2,$ where $P_1, P_2$$\in$ Spec$\left(\mathcal{O}_{K}\right),$ $P_1 \neq P_2$  are non-principal ideals (each of the ideals $P_1$ and $P_2$ has $2$ generators) if and only if $d\equiv 1$ (mod $8$).
In this case $2\mathcal{O}_{K}=P_1 \cdot P_2=P_1 \cap P_2.$ 
\item[$v)$]  $(\mathbf{I(\mathcal{O}_{K}),\subseteq })$ is a Noetherian lattice.
\end{enumerate}
\end{theorem}

\begin{proof} 
i) It follows immediately from Theorem \ref{theorem23.} i) and ii).\\
ii) It follows immediately from Theorem \ref{theorem23.} ii) and from the Chinese Remainder Theorem.\\
iii) It follows immediately from Theorem \ref{theorem23.} iv) and v).\\
iv) It follows immediately from Theorem \ref{theorem23.} v) and from the Chinese Remainder Theorem.\\
v) It is known that $\mathcal{O}_{K}$  is a Dedeking ring and applying Corollary  \ref{corollary36.}, it results that  $(\mathbf{I(\mathcal{O}_{K}),\subseteq })$ is a Noetherian lattice.
\end{proof}

\section{Some unitary rings with all idempotent ideals}

In \cite{ps}, we proved the following result.

\begin{proposition}

\label{prop*} (Proposition 15 from  \cite{ps})  If  $\left(A, +, \cdot \right)$ is a Boolean ring, any ideal of $A$ is idempotent.

\end{proposition}
We wonder if the converse of this Proposition is true and we get that the answer is negative.
We are looking for unitary rings, which are neither fields nor Boolean rings, but which have all idempotent ideals. We found some such examples, which we present below.
\begin{example}  
\label{example 1.}
 $\left(\mathbb{Z}_6, +, \cdot \right)$ is a unitary commutative ring. The lattice of its ideals is $(\mathbf{I(\mathbb{Z}_6),\subseteq }), $ $\mathbf{I(\mathbb{Z}_6)}=\left\{\left(\overline{0}\right), 2\mathbb{Z}/6\mathbb{Z}, 3\mathbb{Z}/6\mathbb{Z},  \mathbb{Z}_6\right\}.$\\
It is clear that $\left(\overline{0}\right)^2 = \left(\overline{0}\right), $ $\mathbb{Z}_6^2=\mathbb{Z}_6.$\\
Let $I=2\mathbb{Z}/6\mathbb{Z}=\left\{\overline{0}, \overline{2}, \overline{4}\right\}.$ We show that $I^2 = I.$ 
$$I^{2}=\left\{x\in \mathbb{Z}_6: x=x_{1}\cdot y_{1}+x_{2}\cdot y_{2}+....+x_{n}\cdot y_{n}|\: n\in \mathbb{N}^{*},\: x_{i},y_{i}\in I, i=1,2,\ldots,n \right\}=$$
$$=\left\{\overline{0}, \overline{0}\cdot \overline{2}, \overline{0}\cdot \overline{4}, \overline{2}\cdot \overline{2}, \overline{2}\cdot \overline{4}, \overline{4}\cdot \overline{4}\right\}=\left\{\overline{0}, \overline{2}, \overline{4}\right\}=I.$$
Let $J=3\mathbb{Z}/6\mathbb{Z}=\left\{\overline{0}, \overline{3}\right\}.$ We show that $J^2 = J.$ 
$$J^{2}=\left\{x\in \mathbb{Z}_6: x= x_{1}\cdot y_{1}+x_{2}\cdot y_{2}+....+x_{n}\cdot y_{n}|\: n\in \mathbb{N}^{*},\: x_{i},y_{i}\in J, i=1,2,\ldots,n \right\}=$$
$$=\left\{\overline{0}, \overline{0}\cdot \overline{3}, \overline{3}\cdot \overline{3}\right\}=\left\{\overline{0}, \overline{3}\right\}=J.$$
So,all ideals of  the ring $\left(\mathbb{Z}_6, +, \cdot \right)$ are idempotent.\\
But $\overline{2}^2 = \overline{4}\neq \overline{2}$ (in $\mathbb{Z}_6$). So,  the ring $\left(\mathbb{Z}_6, +, \cdot \right)$ is not a boolean ring. Since $6$ is not prime, it follows that the ring  $\left(\mathbb{Z}_6, +, \cdot \right)$  is not a field.
\end{example}  
\begin{example}  
\label{example 2.}
 $\left(\mathbb{Z}_{10}, +, \cdot \right)$ is a unitary commutative ring. The lattice of its ideals is $(\mathbf{I(\mathbb{Z}_{10}),\subseteq }), $ $\mathbf{I(\mathbb{Z}_{10})}=\left\{\left(\overline{0}\right), 2\mathbb{Z}/10\mathbb{Z}, 5\mathbb{Z}/10\mathbb{Z},  \mathbb{Z}_{10}\right\}.$\\
It is clear that $\left(\overline{0}\right)^2 = \left(\overline{0}\right), $ $\mathbb{Z}_{10}^2=\mathbb{Z}_{10}.$\\
Let $I=2\mathbb{Z}/10\mathbb{Z}=\left\{\overline{0}, \overline{2}, \overline{4}, \overline{6}, \overline{8} \right\}.$ We prove that $I^2 = I.$ 
$$I^{2}=\left\{x\in \mathbb{Z}_6: x=x_{1}\cdot y_{1}+x_{2}\cdot y_{2}+....+x_{n}\cdot y_{n}|\: n\in \mathbb{N}^{*},\: x_{i},y_{i}\in I, i=1,2,\ldots,n\right\}=$$
$$=\left\{\overline{0}, \overline{0}\cdot \overline{2}, \overline{0}\cdot \overline{4},  \overline{0}\cdot \overline{6}, \overline{0}\cdot \overline{8}, \overline{2}\cdot \overline{2}, \overline{2}\cdot \overline{4},  \overline{2}\cdot \overline{6}, \overline{2}\cdot \overline{8},  \overline{4}\cdot \overline{4},
 \overline{4}\cdot \overline{6},  \overline{4}\cdot \overline{8},  \overline{6}\cdot \overline{6},  \overline{6}\cdot \overline{8},...\right\}=$$
$$=\left\{\overline{0}, \overline{2}, \overline{4}, \overline{6}, \overline{8}\right\}=I.$$
We consider $J=5\mathbb{Z}/10\mathbb{Z}=\left\{\overline{0}, \overline{5}\right\}.$ We prove that $J^2 = J.$ 
$$J^{2}=\left\{x\in \mathbb{Z}_6: x=x_{1}\cdot y_{1}+x_{2}\cdot y_{2}+....+x_{n}\cdot y_{n}|\: n\in \mathbb{N}^{*},\: x_{i},y_{i}\in J, i=1,2,\ldots,n\right\}=$$
$$=\left\{\overline{0}, \overline{0}\cdot \overline{5}, \overline{5}\cdot \overline{5}\right\}=\left\{\overline{0}, \overline{5}\right\}=J.$$
So, all ideals of  the ring $\left(\mathbb{Z}_{10}, +, \cdot \right)$ are idempotent.\\
Since $10$ is not prime, it follows that the ring  $\left(\mathbb{Z}_{10}, +, \cdot \right)$  is not a field.\\
 $\overline{3}^2 = \overline{9}\neq \overline{3}$ (in $\mathbb{Z}_{10}$). So,  the ring $\left(\mathbb{Z}_{10}, +, \cdot \right)$ is not a boolean ring. 

\end{example}  
We managed to generalize the results from these two examples as follows:
\begin{proposition}
\label{prop14.}  Let  $n$ be a positive integer, $n\geq 2,$ $n=p\cdot q$, where $p$ and $q$ are two distinct prime positive integers. Let the unitary commutative ring $\left(\mathbb{Z}_{n}, +, \cdot \right).$ 
 Then:
\begin{enumerate}
\item[$i)$] the ring $\left(\mathbb{Z}_{n}, +, \cdot \right)$ is not a field;

\item[$ii)$]  the ring $\left(\mathbb{Z}_{n}, +, \cdot \right)$ is not a boolean ring;

\item[$ii)$]  all ideals of  the ring $\left(\mathbb{Z}_{n}, +, \cdot \right)$ are idempotent ideals.

\end{enumerate}
\end{proposition}
\begin{proof}
i) Since $n=p\cdot q$ is not a prime number,  it results that the ring  $\left(\mathbb{Z}_{n}, +, \cdot \right)$  is not a field.\\
ii) Since $n=p\cdot q,$  where $p$ and $q$ are two distinct prime positive integers, it results that $n\geq 6$ and from here we get $\overline{2}^{2}\neq \overline{2}$  (in  $\mathbb{Z}_{n}$), so  the ring $\left(\mathbb{Z}_{n}, +, \cdot \right)$ is not a boolean ring. \\
iii) The lattice of ideals of the ring $\left(\mathbb{Z}_{n}, +, \cdot \right)$ is $(\mathbf{I(\mathbb{Z}_{n}),\subseteq }), $ $\mathbf{I(\mathbb{Z}_{n})}=$\\
$\left\{\left(\overline{0}\right), p\mathbb{Z}/n\mathbb{Z}, q\mathbb{Z}/n\mathbb{Z},  \mathbb{Z}_{n}\right\}.$\\
It is clear that $\left(\overline{0}\right)^2 = \left(\overline{0}\right), $ $\mathbb{Z}_{n}^2=\mathbb{Z}_{n}.$\\
We consider the ideals $I=p\mathbb{Z}/n\mathbb{Z},$ $J=q\mathbb{Z}/n\mathbb{Z}.$\\
$I=p\mathbb{Z}/n\mathbb{Z}=\left\{\overline{0}, \overline{p}, \overline{2p}, ..., \overline{\left(q-1\right)\cdot p }\right\}.$ We show that $I^2 = I.$ It is clear that $I^2 \subseteq I.$
We are proving that  $I \subseteq I^2.$ Let $x$$\in$$I,$ $x=\overline{ip},$ where $i\in$$\left\{0,1,...,q-1\right\}.$ We show that there is $k\in$$\left\{0,1,...,q-1\right\}$ such that $x=\overline{kp^2}=\overline{k_{1}p}\cdot \overline{k_{2}p} $$\in$ $I^2,$ 
that is, it is enough to show that there is $k\in$$\left\{0,1,...,q-1\right\}$  such that $i\cdot p$$\equiv$$k\cdot p^2$ (mod $p\cdot q$), so  $k\cdot p$$\equiv$$i$ (mod $q$). This is a congruence in the unknown $k.$
Since $g.c.d.(p,q)=1,$ it results that this congruence has a unique solution modulo $p$ namely $k$$\equiv$$i\cdot p^{\varphi\left(q\right)-1}$ (mod $q$), where $\varphi$ is Euler's function, so $\left(\exists!\right)$ $k\in$$\left\{0,1,...,q-1\right\}$
such that  $x=\overline{k\cdot p^2},$ meaning $x$$\in$$I^2.$ So, $I^2 = I$ Analogously, it shows that $J^2 = J.$ Thus we obtain that  all ideals of  the ring $\left(\mathbb{Z}_{n}, +, \cdot \right)$ are idempotent ideals.
\end{proof}
\smallskip\\
We are interested in whether, for any positive integer $n$, $n\geq 2,$ the ring $\left(\mathbb{Z}_{n}, +, \cdot \right)$ has the property that every ideal of it is indempotent. The answer is negative. We give some examples in this regard.\\

\begin{example} 
\label{example 3.}
 We consider  $n=4$ and the ring $\left(\mathbb{Z}_4, +, \cdot \right).$ The lattice of its ideals is $(\mathbf{I(\mathbb{Z}_4),\subseteq }), $ $\mathbf{I(\mathbb{Z}_4)}=\left\{\left(\overline{0}\right), 2\mathbb{Z}/4\mathbb{Z}, \mathbb{Z}_4\right\}.$\\
It is clear that $\left(\overline{0}\right)^2 = \left(\overline{0}\right), $ $\mathbb{Z}_4^2=\mathbb{Z}_4.$\\
Let the ideal $I=2\mathbb{Z}/4\mathbb{Z}=\left\{\overline{0}, \overline{2}\right\}.$  We have: 
$$I^2= \left\{\displaystyle\sum_{\overline{i}, \overline{j}\in I} \overline{i}\cdot \overline{j} \right\}= \left\{\displaystyle\sum_{\overline{i}, \overline{j}\in \left\{\overline{0}, \overline{2}\right\}} \overline{i}\cdot \overline{j} \right\}=
\left\{\overline{0}\right\}.$$
It results that $I^2\neq I,$ so the ideal $I$ is not idempotent. Thus, the ring  $\left(\mathbb{Z}_4, +, \cdot \right)$ does not have all idempotent ideals.
\end{example} 

\begin{example} 
\label{example 4.}
For  $n=12,$ we consider the ring $\left(\mathbb{Z}_{12}, +, \cdot \right).$ The lattice of its ideals is $(\mathbf{I(\mathbb{Z}_{12}),\subseteq }), $ $\mathbf{I(\mathbb{Z}_{12})}=\left\{\left(\overline{0}\right), 2\mathbb{Z}/12\mathbb{Z}, 3\mathbb{Z}/12\mathbb{Z},  4\mathbb{Z}/12\mathbb{Z},  6\mathbb{Z}/12\mathbb{Z}, \mathbb{Z}_{12}\right\}.$\\
We consider the ideal $J=6\mathbb{Z}/12\mathbb{Z}=\left\{\overline{0}, \overline{6}\right\}.$  We have: 
$$J^2= \left\{\displaystyle\sum_{\overline{i}, \overline{j}\in J} \overline{i}\cdot \overline{j} \right\}= \left\{\displaystyle\sum_{\overline{i}, \overline{j}\in \left\{\overline{0}, \overline{6}\right\}} \overline{i}\cdot \overline{j} \right\}=
\left\{\overline{0}\right\}\neq J.$$
It results that the ideal $J$ is not idempotent. So, the ring  $\left(\mathbb{Z}_{12}, +, \cdot \right)$ does not have all idempotent ideals.
\end{example} 

We obtain the following result.

\begin{proposition}
\label{prop17.}  Let  $n$ be a positive integer, $n\geq 2,$ $n=p^2\cdot q$, where $p$ and $q$ are two distinct prime positive integers. Let the unitary commutative ring $\left(\mathbb{Z}_{n}, +, \cdot \right).$ 
 Then:
\begin{enumerate}
\item[$i)$] the ring $\left(\mathbb{Z}_{n}, +, \cdot \right)$ is  not a field;

\item[$ii)$]  the ring $\left(\mathbb{Z}_{n}, +, \cdot \right)$ is  not a boolean ring;

\item[$iii)$] The ring  $\left(\mathbb{Z}_{n}, +, \cdot \right)$  does not have all idempotent ideals.

\end{enumerate}
\end{proposition}
\begin{proof}
i) Since $n=p^2\cdot q$ is not a prime number,  it results that the ring  $\left(\mathbb{Z}_{n}, +, \cdot \right)$  is not a field.\\
ii) Since $n=p^2\cdot q,$  where $p$ and $q$ are two distinct prime positive integers, it results that $n\geq 12$ and from here we get $\overline{2}^{2}\neq \overline{2}$  (in  $\mathbb{Z}_{n}$), so  the ring $\left(\mathbb{Z}_{n}, +, \cdot \right)$ is not not a boolean ring. \\
iii) In the ring $\left(\mathbb{Z}_{n}, +, \cdot \right)$ we consider the ideal  
$$I=p\cdot q\mathbb{Z}/n\mathbb{Z}=  p\cdot q\mathbb{Z}/ p^2\cdot q\mathbb{Z}= \left\{\overline{0}, \overline{p\cdot q}, \overline{2\cdot p\cdot q},..., \overline{\left(p-1\right)\cdot p\cdot q} \right\}.$$
It results that
$$I^2= \left\{\displaystyle\sum_{\overline{i}, \overline{j}\in \left\{\overline{0}, \overline{1},..., \overline{p-1} \right\}} \overline{i\cdot  p\cdot q}\cdot \overline{j\cdot  p\cdot q} \right\}= \left\{{\overline{0}} \right\} \neq I,$$
so the ideal $I$ is not idempotent. We obtain that the ring $\left(\mathbb{Z}_{n}, +, \cdot \right)$  does not have all idempotent ideals.
\end{proof}
\smallskip\\
We wonder for which natural numbers $n,$  $n\geq 2,$  the ring $\left(\mathbb{Z}_{n}, +, \cdot \right)$ simultaneously has the following properties: it has all idempotent ideals, it is not a field and it is not a boolean ring.
We obtain the following result.

\begin{proposition}
\label{prop18.}  Let  $n$ be a positive integer, $n\geq 4,$ $n$ is not a prime. 
 Then:
\begin{enumerate}
\item[$i)$] the ring $\left(\mathbb{Z}_{n}, +, \cdot \right)$ is not a field;

\item[$ii)$]  the ring $\left(\mathbb{Z}_{n}, +, \cdot \right)$ is not a boolean ring;

\item[$iii)$] The ring  $\left(\mathbb{Z}_{n}, +, \cdot \right)$  has all idempotent ideals if and only if $n$ is square free.

\end{enumerate}
\end{proposition}
\begin{proof}
i) Since $n$ is not a prime number,  it results that the ring  $\left(\mathbb{Z}_{n}, +, \cdot \right)$  is not a field.\\
ii) Since  $n\geq 4,$ , it results that $\overline{2}^{2}\neq \overline{2}$  (in  $\mathbb{Z}_{n}$), so  the ring $\left(\mathbb{Z}_{n}, +, \cdot \right)$ is not a boolean ring. \\
iii) $"\Rightarrow"$ We assume that the ring  $\left(\mathbb{Z}_{n}, +, \cdot \right)$  has all idempotent ideals and we prove that $n$ is a square free integer.\\
Let assume by absurd that there exists $n$$\in$$\mathbb{N}$,  $n\geq 4,$ which is not prime and $n$ is not square free, such that all ideals of the ring $\left(\mathbb{Z}_{n}, +, \cdot \right)$ are idempotent.\\
Since $n$ is not square free it results that there is $p_1$ a prime positive integer such that $p^{2}_{1} | n.$ According to the Fundamental Theorem of Arithmetic, $n$ is written uniquely:
$n= p^{\alpha_{1}}_{1}\cdot p^{\alpha_{2}}_{2}\cdot...\cdot p^{\alpha_{r}}_{r},$ where $p_1,$ $p_2,$..., $p_r$ are distinct
prime integers and $\alpha_{1}, \alpha_{2},..., \alpha_{r}$$\in$$\mathbb{N}^{*},$ $\alpha_{1}\geq 2.$\\
We consider the following ideals of the ring $\left(\mathbb{Z}_{n}, +, \cdot \right):$ 
$$I=p^{\alpha_{1}-1}_{1}\cdot p^{\alpha_{2}}_{2}\cdot...\cdot p^{\alpha_{r}}_{r}\cdot \mathbb{Z}/n\cdot \mathbb{Z}= p^{\alpha_{1}-1}_{1}\cdot p^{\alpha_{2}}_{2}\cdot...\cdot p^{\alpha_{r}}_{r}\cdot \mathbb{Z}/ p^{\alpha_{1}}_{1}\cdot p^{\alpha_{2}}_{2}\cdot...\cdot p^{\alpha_{r}}_{r}\cdot \mathbb{Z}=$$
$$=\left\{\overline{0}, \overline{p^{\alpha_{1}-1}_{1}\cdot p^{\alpha_{2}}_{2}\cdot...\cdot p^{\alpha_{r}}_{r}},  \overline{2\cdot p^{\alpha_{1}-1}_{1}\cdot p^{\alpha_{2}}_{2}\cdot...\cdot p^{\alpha_{r}}_{r}},  ... ,  \overline{\left(p_{1} -1\right)\cdot p^{\alpha_{1}-1}_{1}\cdot p^{\alpha_{2}}_{2}\cdot...\cdot p^{\alpha_{r}}_{r}}  \right\}.$$
It results that
$$I^2= \left\{\displaystyle\sum_{\overline{a}, \overline{b}\in I} \overline{a}\cdot \overline{b} \right\}= \left\{\displaystyle\sum_{\overline{i}, \overline{j}\in \left\{\overline{0}, \overline{1},..., \overline{p_{1}-1} \right\}} \overline{i\cdot p^{\alpha_{1}-1}_{1}\cdot p^{\alpha_{2}}_{2}\cdot...\cdot p^{\alpha_{r}}_{r}}\cdot \overline{j\cdot p^{\alpha_{1}-1}_{1}\cdot p^{\alpha_{2}}_{2}\cdot...\cdot p^{\alpha_{r}}_{r}} \right\}=$$
$$=\left\{\displaystyle\sum_{\overline{i}, \overline{j}\in \left\{\overline{0}, \overline{1},..., \overline{p_{1}-1} \right\}} \overline{i}\cdot  \overline{j}\cdot  \overline{p^{2\alpha_{1}-2}_{1}\cdot p^{2\alpha_{2}}_{2}\cdot...\cdot p^{2\alpha_{r}}_{r}} \right\},$$
Since $\alpha_{1}\geq 2,$ it results that  $2\alpha_{1}-2\geq \alpha_{1},$ so $ \overline{p^{2\alpha_{1}-2}_{1}\cdot p^{2\alpha_{2}}_{2}\cdot...\cdot p^{2\alpha_{r}}_{r}}=  \overline{0} $ in $\mathbb{Z}_{n}.$ So $I^2= \left(\overline{0}\right) \neq I.$
We obtain that the ideal $I$ is not idempotent. It follows that the assumption we made is false, so $n$ is square free.\\

$"\Leftarrow"$ Let  $n$ be a positive integer, $n\geq 4,$ $n$ is not a prime, $n$ is square free. It results that $n=p_{1}\cdot p_{2}\cdot...\cdot p_{l},$  where $l\in \mathbb{N},$ $l\geq 2$ and $p_{1}, p_{2},...,p_{l}$ are prime distinct positive integers. 
The lattice of ideals of the ring $\left(\mathbb{Z}_{n}, +, \cdot \right)$ is $(\mathbf{I(\mathbb{Z}_{n}),\subseteq }), $ 

$$\mathbf{I(\mathbb{Z}_{n})}=\left\{ d\mathbb{Z}/p_{1}\cdot p_{2}\cdot...\cdot p_{l}\mathbb{Z}, \  where  \  d | p_{1}\cdot p_{2}\cdot...\cdot p_{l} \right\}=$$
$$=\left\{ \left(\overline{0}\right),  p_{1}\mathbb{Z}/p_{1}\cdot p_{2}\cdot...\cdot p_{l}\mathbb{Z}, ..., p_{2}\cdot... \cdot.p_{l}\mathbb{Z}/p_{1}\cdot p_{2}\cdot...\cdot p_{l}\mathbb{Z} \right\}.$$
Since $n=p_{1}\cdot p_{2}\cdot...\cdot p_{l},$ it results that the number of natural divisors of $n$ is $\tau\left(n\right)=2^{l},$ so Card$\mathbf{I(\mathbb{Z}_{n})}=2^{l}.$\\
Let $I$$\in$$\mathbf{I(\mathbb{Z}_{n})}.$ It results that $I=p_{i_{1}}\cdot... \cdot.p_{i_{s}}\mathbb{Z}/p_{1}\cdot p_{2}\cdot...\cdot p_{l}\mathbb{Z},$ where $s\in \mathbb{N}^{*},$ $s\leq l,$ $ \left\{p_{i_{1}},...,p_{i_{s}}\right\}\subseteq $$ \left\{p_{1},...,p_{l}\right\}. $
Without reducing the generality, we can rename the indices so that  $ \left\{p_{i_{1}},...,p_{i_{s}}\right\}=\left\{p_{1},...,p_{s}\right\}.$ It results that the ideal $I=p_{1}\cdot... \cdot.p_{s}\mathbb{Z}/p_{1}\cdot p_{2}\cdot...\cdot p_{l}\mathbb{Z},$  where $s\in \mathbb{N}^{*},$ $s\leq l.$\\
We prove that $I^2 = I.$ It is clear that $I^2 \subseteq I.$ We are showing that $I \subseteq I^2.$ Let
 $$x\in I=p_{1}\cdot... \cdot.p_{s}\mathbb{Z}/p_{1}\cdot p_{2}\cdot...\cdot p_{l}\mathbb{Z}=$$
$$=\left\{\overline{0}, \overline{p_{1}\cdot p_{2}\cdot...\cdot p_{s}}, \overline{2\cdot p_{1}\cdot p_{2}\cdot...\cdot p_{s}},
..., \overline{\left(p_{s+1}\cdot...\cdot p_{l}-1\right)\cdot p_{1}\cdot...\cdot p_{s}}\right\}.$$
It results that $x= \overline{j\cdot p_{1}\cdot p_{2}\cdot...\cdot p_{s}},$ where $j$$\in$$\left\{0, 1, 2,..., p_{s+1}\cdot...\cdot p_{l}-1\right\}.$\\
 Now, we are prove that there is $k\in$$\left\{0, 1, 2,..., p_{s+1}\cdot...\cdot p_{l}-1\right\}$ such that 
$$x=\overline{k\cdot p^{2}_{1}\cdot p^{2}_{2}\cdot...\cdot p^{2}_{s}}=\overline{k_{1}\cdot p_{1}\cdot p_{2}\cdot...\cdot p_{s}}\cdot
 \overline{k_{2}\cdot p_{1}\cdot p_{2}\cdot...\cdot p_{s}}$$
(meaning $x$ is in $I^2$).\\
So, it is enough to show that there is $k\in$$\left\{0, 1, 2,..., p_{s+1}\cdot...\cdot p_{l}-1\right\}$ such that $ \overline{j\cdot p_{1}\cdot p_{2}\cdot...\cdot p_{s}}=\overline{k\cdot p^{2}_{1}\cdot p^{2}_{2}\cdot...\cdot p^{2}_{s}}$ (in $\mathbb{Z}_{n}$), which is equivalent to the congruence 
$$j\cdot p_{1}\cdot p_{2}\cdot...\cdot p_{s} \equiv k\cdot p^{2}_{1}\cdot p^{2}_{2}\cdot...\cdot p^{2}_{s}\ (mod \ p_{1}\cdot p_{1}\cdot...\cdot p_{s}\cdot...\cdot p_{l})  $$
has solutions, i.e., the following congruence (in the unknown $k$)
$$k\cdot p_{1}\cdot p_{2}\cdot...\cdot p_{s} \equiv j\ (mod \ p_{s+1}\cdot p_{s+2}\cdot...\cdot p_{l})  $$
has solutions. This is true, from theory, because\\
 gcd$\left(p_{1}\cdot p_{2}\cdot...\cdot p_{s}, \ p_{s+1}\cdot p_{s+2}\cdot...\cdot p_{l}\right)=1.$ It results $x\in I^2 .$ So  $I \subseteq I^2.$
\end{proof}

\begin{proposition}
\label{prop19.}  Let  $l$ be a positive integer, $l\geq 2$ and let $M$ be a finite set with $l$ elements, $M=\left\{p_{1}, p_{2},..., p_{l}\right\},$ where $p_{1}, p_{2},..., p_{l}$ are distinct prime positive integers. Then, the lattice 
$\left(P\left(M\right), \subseteq, \vee, \wedge\right)$ is isomorphic with the lattice of ideals of the ring $\left(\mathbb{Z}_{p_{1}\cdot p_{2}\cdot ...\cdot  p_{l}}, +, \cdot \right).$ 
\end{proposition}
\begin{proof} Let $f:P\left(M\right)\longrightarrow$$\mathbf{I(}\left(\mathbb{Z}_{p_{1}\cdot p_{2}\cdot ...\cdot  p_{l}}, +, \cdot \right)),$\\
 $f\left(\left\{p_{i}\right\} \right)=p_{1}\cdot p_{2}\cdot...\cdot p_{i-1}\cdot  p_{i+1}\cdot...\cdot p_{l}\cdot \mathbb{Z}/p_{1}\cdot p_{2}\cdot...\cdot p_{l}\cdot \mathbb{Z},$ $\forall$ $i=$$1,\ldots, l,$\\
$f\left(\left\{p_{i_{1}}, p_{i_{2}},..., p_{i_{r}}\right\}\right)=p_{i_{r+1}}\cdot p_{i_{r+2}}\cdot...\cdot p_{l}\cdot \mathbb{Z}/p_{1}\cdot p_{2}\cdot...\cdot p_{l}\cdot \mathbb{Z},$ $\forall$ $r\in \mathbb{N}*$, $r\leq l,$\\
$\forall$ $\left\{p_{i_{1}}, p_{i_{2}},..., p_{i_{r}}\right\}$$\subseteq$$\left\{p_{1}, p_{2},..., p_{l}\right\}.$\\
If $p_{i}, p_{j}$$\in$$\left\{p_{1}, p_{2},..., p_{l}\right\},$ $p_{i}\neq p_{j}, $ it results that the ideals
 $$p_{1}\cdot p_{2}\cdot...\cdot p_{i-1}\cdot  p_{i+1}\cdot...\cdot p_{l}\cdot \mathbb{Z}/p_{1}\cdot p_{2}\cdot...\cdot p_{l}\cdot \mathbb{Z}\neq p_{1}\cdot p_{2}\cdot...\cdot p_{j-1}\cdot  p_{j+1}\cdot...\cdot p_{l}\cdot \mathbb{Z}/p_{1}\cdot p_{2}\cdot...\cdot p_{l}\cdot \mathbb{Z}$$
$$ \Leftrightarrow  f\left(\left\{p_{i}\right\} \right) \neq  f\left(\left\{p_{j}\right\} \right).$$
Analogously, we get that $\forall$ $ M_{1}, M_{2}\subseteq M, $ $M_{1}\neq M_{2},$ it results that $f\left(M_{1}\right)\neq f\left(M_{2}\right).$ So, $f$ is injective. Moreover Card$P\left(M\right)=2^{l}$$=Card \mathbf{I(}\left(\mathbb{Z}_{p_{1}\cdot p_{2}\cdot ...\cdot  p_{l}}, +, \cdot \right)),$ 
It results that $f$ is bijective.\\
We are showing that $f$ is a lattice morphism.\\
Let $A, B$$\in$$P\left(M\right),$ $A=\left\{p_{i_{1}}, p_{i_{2}},..., p_{i_{r}}\right\}\subseteq$$\left\{p_{i_{1}}, p_{i_{2}},..., p_{i_{r}}, p_{i_{r+1}},...p_{i_{s}}\right\}=B. $ It results that $f\left(A\right)=$$ p_{i_{r+1}}\cdot...\cdot p_{s}\cdot  p_{s+1}\cdot...\cdot p_{l}\cdot \mathbb{Z}/p_{1}\cdot p_{2}\cdot...\cdot p_{l}\cdot \mathbb{Z}$ and $f\left(B\right)=$$p_{s+1}\cdot...\cdot p_{l}\cdot \mathbb{Z}/p_{1}\cdot p_{2}\cdot...\cdot p_{l}\cdot \mathbb{Z}.$ Since $p_{s+1}\cdot...\cdot p_{l} | p_{i_{r+1}}\cdot...\cdot p_{s}\cdot  p_{s+1}\cdot...\cdot p_{l},$ it results that
$ p_{i_{r+1}}\cdot...\cdot p_{s}\cdot  p_{s+1}\cdot...\cdot p_{l}\cdot \mathbb{Z} $$\subseteq$ $p_{s+1}\cdot...\cdot p_{l}\mathbb{Z} $ and so $f\left(A\right) \subseteq f\left(B\right).$ It follows that  $f$ is a lattice morphism.\\
Since $f$ is bijective and $f$ is a lattice morphism, we obtain that $f$ is a lattice isomorphism.
\end{proof}
\smallskip\\
Finally, using Proposition \ref{prop19.} we can give another proof of $"\Leftarrow"$  of Proposition  \ref{prop18.}, namely: according to Proposition \ref{prop19.}, the lattice  $\left(P\left(M\right), \subseteq, \vee, \wedge\right)$ (where  $M=\left\{p_{1}, p_{2},..., p_{l}\right\}$) is isomorphic to the lattice of ideals of the ring     $\left(\mathbb{Z}_{p_{1}\cdot p_{2}\cdot ...\cdot  p_{l}}, +, \cdot \right).$  But, as every element of the lattice $\left(P\left(M\right), \subseteq, \vee, \wedge\right)$ is idempotent, it follows that every ideal of the ring $\left(\mathbb{Z}_{p_{1}\cdot p_{2}\cdot ...\cdot  p_{l}}, +, \cdot \right)$ is idempotent.\\
\smallskip\\
We wonder if there are rings other than $\left(\mathbb{Z}_{n}, +, \cdot \right)$ which are not fields, not boolean rings but their ideals are idempotent. The answer is affirmative. A simple example of a ring with these properties is the commutative ring  $\left(\mathbb{Z}_{p}\times \mathbb{Z}_{q}, +, \cdot \right)$ (which is isomorphic to the ring $\left(\mathbb{Z}_{pq}, +, \cdot \right)$), with p and q prime positive integers, $p\neq q.$\\
In the following we give other examples of such rings.\\
The following remark is known (see \cite{DF}).

\begin{remark}  
\label{remark*.}   Let $A$ be a commutative and unitary ring. 
\begin{enumerate}
\item[$i)$] If $I$ is an ideal of $A$, then $M_{n}\left(I\right)$ is a two-sided ideal of the  non-commutative and unitary ring $M_{n}\left(A\right)$;
\item[$ii)$] If $J$ is a two-sided ideal of the  non-commutative and unitary ring $M_{n}\left(A\right)$ then there is some ideal $I$ of $A$
such that $J =M_{n}\left(I\right)$.
\end{enumerate}
\end{remark}
Using this Remark and Proposition \ref{prop18.} we obtain the following result.
\begin{proposition}
\label{prop20.} Let  $n, m$ be two positive integers, $n\geq 4,$ $n$ is not a prime, $m\geq 2$. Let  $\left( M_{m}\left(\mathbb{Z}_{n}\right),  +, \cdot \right)$ be the non-commutative and unitary ring of square matrices of order $m$ with elements from $\mathbb{Z}_{n}.$
 Then:
\begin{enumerate}
\item[$i)$] the ring $\left( M_{m}\left(\mathbb{Z}_{n}\right),  +, \cdot \right)$ is not a field;

\item[$ii)$]  the ring $\left( M_{m}\left(\mathbb{Z}_{n}\right),  +, \cdot \right)$ is not a boolean ring;

\item[$iii)$] the ring $\left( M_{m}\left(\mathbb{Z}_{n}\right),  +, \cdot \right)$ has all idempotent two-sided ideals if and only if $n$ is square free.

\end{enumerate}
\end{proposition}
\begin{proof} i) Let the matrix $A=\left( 
\begin{array}{ccccc}
\overline{1} & \overline{1} & \overline{1} & ... & \overline{1} \\ 
\overline{0} & \overline{0} & \overline{0} & ... & \overline{0} \\ 
\overline{0} & \overline{0} & \overline{0} & ... & \overline{0} \\ 
... & ... & ... & ... & ... \\ 
\overline{0} & \overline{0} & ... & \overline{0} & \overline{0}%
\end{array}%
\right)  $$\in$$M_{m}\left(\mathbb{Z}_{n}\right).$ Since det$\left(A\right)=\overline{0},$ it follows that $A$ is not invertible, so the ring $M_{m}\left(\mathbb{Z}_{n}\right)$ is not a field.\\
ii) Let the matrix $B=\left( 
\begin{array}{ccccc}
\overline{1} & \overline{1} & \overline{0} & ... & \overline{0} \\ 
\overline{1} & \overline{0} & \overline{0} & ... & \overline{0} \\ 
\overline{0} & \overline{0} & \overline{0} & ... & \overline{0} \\ 
... & ... & ... & ... & ... \\ 
\overline{0} & \overline{0} & ... & \overline{0} & \overline{0}%
\end{array}%
\right)  $$\in$$M_{m}\left(\mathbb{Z}_{n}\right).$\\
Since $B^{2}=\left( 
\begin{array}{ccccc}
\overline{2} & \overline{1} & \overline{0} & ... & \overline{0} \\ 
\overline{1} & \overline{1} & \overline{0} & ... & \overline{0} \\ 
\overline{0} & \overline{0} & \overline{0} & ... & \overline{0} \\ 
... & ... & ... & ... & ... \\ 
\overline{0} & \overline{0} & ... & \overline{0} & \overline{0}%
\end{array}%
\right)$ $\neq B,$ it follows that $M_{m}\left(\mathbb{Z}_{n}\right)$ is not a boolean ring.\\
iii) Let let $\mathbf{I(\mathbb{Z}_{n})}$ be the lattice of ideals of the ring $\left(\mathbb{Z}_{n}, +, \cdot \right)$ and let $\mathbf{I(M_{m}\left(\mathbb{Z}_{n}\right))}$ be the lattice of two-sided ideals of the ring $\left( M_{m}\left(\mathbb{Z}_{n}\right),  +, \cdot \right)$.\\
Let the function $f:$$\mathbf{I(\mathbb{Z}_{n})}$$\rightarrow$$\mathbf{I(M_{m}\left(\mathbb{Z}_{n}\right))}$, $f\left(I\right)=$$M_{m}\left(I\right),$ for each ideal $I$ of the ring $\left(\mathbb{Z}_{n}, +, \cdot \right)$. It immediately follows that $f$ is injective. Applying Remark \ref{remark*.} we obtain that $f$ is bijective.\\
Let $I$ and $J$ be two ideals of the ring $\mathbb{Z}_{n}$. We have that $M_{m}\left(I\right)$  and $M_{m}\left(J\right)$ are two-sided ideals  ideals of the ring $M_{m}\left(\mathbb{Z}_{n}\right)$ and
$$M_{m}\left(I\right)\cdot M_{m}\left(J\right)=$$
$$= \left\{A_{1}\cdot B_{1}+\ldots A_{k}\cdot B_{k} |k\in \mathbb{N}^{*}, A_{1},\ldots, A_{k}\in M_{m}\left(I\right), B_{1},\ldots, B_{k}\in M_{m}\left(J\right)
\right\}$$
$$=M_{m}\left(I\cdot J\right).$$
This implies that $\left(M_{m}\left(I\right)\right)^{2}=M_{m}\left(I^{2}\right),$ for any ideal $I$ of the ring $\left(\mathbb{Z}_{n},  +, \cdot \right).$ Using this equality, the fact that $f$ is bijective and applying Proposition  \ref{prop18.}  we obtain that the ring $\left( M_{m}\left(\mathbb{Z}_{n}\right),  +, \cdot \right)$   has all idempotent two-sided ideals if and only if $n$ is square free.
\end{proof} \\
\smallskip\\
\textbf{Conflicts of Interest}: The author declare no conflict of interest.\\
\smallskip\\
\textbf{Acknowledgements.}
The author wishs express her gratitude to the anonymous reviewers for their careful reading of
the manuscript and for many valuable suggestions.
\begin{equation*}
\end{equation*}

\end{document}